\documentclass[12pt,fleqn]{amsart}
\usepackage{amsmath}
\usepackage{amssymb}
\usepackage{amsthm}
\usepackage[margin=1in]{geometry}
\usepackage{hhline}
\usepackage{setspace}

\usepackage{graphicx}
\usepackage{epstopdf}
\usepackage{xcolor}

\usepackage{geometry}

\newcommand{\ul}{\pmb}
\def\ulul#1{\mathbf{#1}}

\title{Time evolution of a mean-field generalized contact process}
\author{by Logan Chariker$^1$ and Joel L. Lebowitz$^2$}
\thanks{$^1$School of Natural Sciences, Institute for Advanced Study, Princeton, New Jersey. \newline $^2$Departments of Mathematics and Physics, Rutgers University, New Brunswick, New Jersey. \newline\newline Keywords: generalized contact process, neurons with discrete voltage, mean field, exact solution, spatial dependence}

\begin{document}

\maketitle

\noindent{\bf Abstract.}
   We investigate the macroscopic time evolution and stationary states of a mean field generalized contact process in $\mathbb{R}^d$. The model is described by a coupled set of nonlinear integral-differential equations.  It was inspired by a model of neurons with discrete voltages evolving by a stochastic integrate and fire mechanism.  We obtain a complete solution in the spatially uniform case
and partial solutions in the general case.  The system has one or more fixed points and also traveling wave solutions.

\section{Introduction}
We consider the mean-field (hydrodynamic) limit of a (novel) stochastic lattice system inspired by neuronal integrate-and-fire models \cite{Gerstner}.  On the microscopic level the system consists of variables $S(\ul{z},t)$ associated to lattice site $\ul{z}\in\Omega\subset\mathbb{Z}^d$, at time $t\in\mathbb{R}_+$.  $S(\ul{z},t)$ can take integer values $0,1,\ldots,k$. These can be thought of as discrete values of the voltage of a neuron at site $\ul{z}$, or as the state of infection of the individual at $\ul{z}$. There are $k$ stages of the infection with $S(\ul{z},t)=0$ corresponding to the healthy state.

A sketch of the microscopic dynamics is as follows: when none of the $S(\ul{z},t)$ are in state $k$, the system is in a static state.  When $S(\ul{z},t)=k$ it has a probability $dt$ of firing (healing, dying) during the time interval $(t, t+dt)$.  This is independent of the values $S(\ul{w},t)$ for $\ul{w}\neq \ul{z}$.  When $S(\ul{z},t)$ fires it is instantaneously reset to $S'(\ul{z},t)=0$.  A neuron at site $\ul{w}$, with $S(\ul{w},t)=j$, $j\leq k-1$, will jump from $j$ to $j+1$ with rate $k\lambda J_\gamma(\ul{z},\ul{w})$, $\lambda > 0$, when $S(\ul{z},t)=k$.  The function $J_\gamma (\ul{x},\ul{y})$ has the Kac form \cite{Giacomin},
\begin{align}
J_\gamma (\ul{z},\ul{w}) = \gamma^d J(\gamma (\ul{z}-\ul{w})) = \gamma^d J(\gamma (\ul{w}-\ul{z})) \geq 0 \label{eq:Kac}
\end{align}
with $\int_{\mathbb{R}^d} J(\ul{r}) d\ul{r} = 1$.

When $k=1$, each site can exist in one of only two states, one an inactive state and the other active, and in this respect the model becomes similar to a contact process and also to a popular two-state neuron system known as the stochastic Wilson-Cowan model \cite{Goychuk,Zankoc}.  The case $k>1$ introduces multiple inactive states $0,\ldots, k-1$, in which a site exerts no influence on its neighbors, and it must traverse sequentially through the states to reach the active state $k$, which is followed by a reset at rate 1 to the state $0$.  This is inspired by integrate and fire neuron models, where the firing of other neurons is required to drive the membrane potential of a particular neuron from a resting potential to a threshold level in order for it to become active itself.  There is an abundance of integrate and fire neuron models used throughout neuroscience.  Less tractable than two-state neuron models, many results for integrate and fire models are given by numerical simulation (e.g., \cite{Milton}, \cite{Chariker}, \cite{Protachevicz}), although there are some results proved, in particular for non-spatially-dependent networks in the mean field limit\cite{Brunel}.

We are interested here in the inclusion of a spatial structure. \cite{DeMasiLocherbach} and \cite{Duarte} have shown rigorous results in the hydrodynamic limit of a network of neurons, which have continuous membrane potentials and a stochastic firing threshold.  Individual neuron firings produce infinitesimal jumps in the membrane potential of all other neurons, weighted by a spatially dependent connectivity function.  In contrast, neurons in the model studied in this paper have a fixed firing threshold.

We will be interested in the macroscopic equations on the scale $x=\gamma z$ in the mean field limit $\gamma \rightarrow 0$.  We shall not discuss here the derivation of the macroscopic equations from the microscopic models.  That will be done in a different publication \cite{DeMasi}.  Here we shall discuss the solution of the resulting macroscopic equations for different values of $k$.

\section{Macroscopic Equations}

The $\gamma\rightarrow 0$ limit of the microscopic model yields the following equations on the macroscopic spatial scale for the $v_j(\ul{x},t)$, the fraction of the population density at position $x$ in state $j$, $j=\{0,1,\ldots,k\}$,
\begin{align}
\frac{\partial v_0(\ul{x},t)}{\partial t} &= v_k(\ul{x},t) -v_0(\ul{x},t) \lambda k R_k(\ul{x},t)  \label{eq:4} \\
\frac{\partial v_j(\ul{x},t)}{\partial t} &= [v_{j-1}(\ul{x},t) - v_j(\ul{x},t)] \lambda k R_k(\ul{x},t), \quad j=1,\ldots,k-1 \label{eq:5} \\
\frac{\partial v_k(\ul{x},t)}{\partial t} &= -v_k(\ul{x},t) +v_{k-1}(\ul{x},t)\lambda k R_k(\ul{x},t), \label{eq:6}
\end{align}
where $\ul{x}\in\Lambda\subset\mathbb{R}^d$ is a cubical box of sides $L$, with periodic boundary conditions and uniform density equal to one,
\begin{align}
R_k(\ul{x},t) = \int_{\Lambda}J(\ul{x}-\ul{y})v_k(\ul{y},t)d\ul{y}. \label{eq:Rk}
\end{align}
We assume that $J(r)$ has a range less than $L/2$. In the spatially uniform state, when $v_j(\ul{y},t)$ is independent of $\ul{y}$, then, by (\ref{eq:Rk}) and (\ref{eq:Kac}), $R_k=v_k$.

It follows from equations (\ref{eq:4})-(\ref{eq:6}) that starting with $v_j(\ul{x},0)\geq 0$, $\sum_{j=0}^k v_j(\ul{x},0)=1$, then
\begin{align}
\sum_{i=0}^k v_i(\ul{x},t) = \sum_{i=0}^k v_i(\ul{x},0) = 1, \quad v_i(\ul{x},t)\geq 0, \qquad \text{for all $t\geq 0$}. \label{eq:6a}
\end{align}

It is clear from the above equations that if $v_k(\ul{x},0)=0$ for all $\ul{x}$ then $v_i(\ul{x},t)=v_i(\ul{x},0)$ and the system remains in its initial state forever.  This is not very interesting and we shall assume from now on that $v_k(\ul{x},0)>0$ for some values of $\ul{x}$.  It is still possible however that $v_k(\ul{x},t)\rightarrow 0$ as $t\rightarrow \infty$.  We shall call that ``extinction''.  

Using (\ref{eq:6a}) we can replace $v_k(x,t)$ in (2) and (3) by 
\[
v_k(x,t)=1-\sum_{j=0}^{k-1} v_j(x,t)
\]
leading to a closed set of equations for $v_j(x,t)$, $j=0,\ldots,k-1$ with 
\[
R_k(x,t) = \int_\Lambda J(x-y)\Big[1 -\sum_{j=0}^{k-1} v_j(y,t)\Big]dy.
\]

\section{Stationary States}

Consider now the stationary solutions of (\ref{eq:4})-(\ref{eq:6}).  We see that in addition to the $R_k=0$ solution corresponding to $v_k=0$, there is a stationary solution $R_k(\ul{x},t) = \bar{R}(\ul{x}) > 0$ of the form 
\begin{align}
\bar{v}_j (\ul{x}) = \frac{1}{\lambda k \bar{R}(\ul{x}) + k}, \quad j\leq k-1 \label{eq:10p}
\end{align}
while
\begin{align}
\bar{v}_k (\ul{x}) &= \frac{\lambda \bar{R}(\ul{x})}{1 + \lambda \bar{R}(\ul{x})}, \label{eq:11p}
\end{align}
and $\bar{R}(\ul{x})$ satisfies the eq.
\begin{align}
\bar{R}(\ul{x}) &= \lambda\int J(\ul{x}-\ul{y})\frac{ \bar{R}(\ul{y}) }{ \lambda\bar{R}(\ul{y})+1 } d\ul{y}. \label{eq:9b}
\end{align}

Equations (\ref{eq:10p})-(\ref{eq:9b}) have the spatially uniform solution
\begin{align}
\bar{R}=\bar{v}_k = \frac{\lambda -1}{\lambda} \quad \text{ and }\quad \bar{v}_j = \frac{1}{\lambda k}, \quad j<k,\label{eq:13p}
\end{align}
which is a physical sustaining stationary solution with $\bar{v}_k>0$, for $\lambda > 1$. For $\lambda \leq 1$ the only physical stationary solution is $\bar{v}_k=0$. 

Note that equation (\ref{eq:9b}) for the stationary $\bar{R}(\ul{x})$ is independent of $k$.  Thus if there exists a non-vanishing, spatially dependent stationary state $\bar{R}(\ul{x})$, then it will be so for all $k$.  We will show in section~\ref{sec:6} that for $k=1$, the only stationary nonzero $v_k$ is the spatially uniform one.  Hence this will be true for all $k$.

The stationary state $\bar{\ul{v}} = (\bar{v}_0,\ldots,\bar{v}_k)$ given in (\ref{eq:13p}) is linearly stable (as shown in section \ref{sec:5}), so if we start close enough to $\bar{\ul{v}}$ then the system will always approach $\bar{\ul{v}}$ as $t\rightarrow \infty$.  On the other hand there are, as shown in section \ref{sec:4} for $k>1$, initial uniform states $v_j(0)$ with $v_k(0)>0$ such that $v_k(t)\rightarrow 0$ as $t\rightarrow\infty$.  Initial states $v_j(\ul{x},0)$ close to $v_j(0)$ would also have $v_k(\ul{x},t)\rightarrow 0$ as $t\rightarrow\infty$ for all $\ul{x}$; see figure \ref{fig:3} and section~\ref{sec:concl}.

\section{Solution of the macroscopic equations in the spatially uniform case}\label{sec:4}

The macroscopic equations in the spatially uniform case take the form
\begin{align}
  \frac{dv_0}{dt} &= v_k - (k\lambda v_k)v_0 \label{eq:1} \\
  \frac{dv_j}{dt} &= (k\lambda v_k)[v_{j-1}-v_j], \qquad j=1,\ldots,k-1 \\
  \frac{dv_k}{dt} &= -v_k + (k\lambda v_k)v_{k-1} \label{eq:3}
\end{align}
where $v_j(t)=$ fraction of neurons in state $j=0,1,\ldots,k$.

To simplify equations (\ref{eq:1})-(\ref{eq:3}), we introduce the variable $r(t)$ defined by the equations
\[
  \frac{dr}{dt}=k\lambda v_k,\qquad r(0)=0.
\]
Writing $v_j(t) = \tilde{v}_j(r(t))$, valid as long as $v_k>0$, we get
\begin{align}
  \frac{d\tilde{v}_0}{dr} &= -\tilde{v}_0 + q_0, \qquad q_0 = \frac{1}{k\lambda}, \\
  \frac{d\tilde{v}_j}{dr} &= \tilde{v}_{j-1} - \tilde{v}_j, \qquad j=1,\ldots,k-1, \\
  \frac{d\tilde{v}_k}{dr} &= \tilde{v}_{k-1} + q_k, \qquad q_k=-q_0=-\frac{1}{\lambda k}.
\end{align}  

The autonomous equations
 for $j=0,1,\ldots,k-1$ can be rewritten in vector form
\begin{align}
  \frac{d\tilde{\ul{v}}}{dr} &= \ulul{A}\tilde{\ul{v}} + \ul{q}, \qquad \tilde{\ul{v}}(0) = \ul{v}(0) \label{eq:vec}
\end{align}
where $\tilde{\ul{v}} = (\tilde{v}_0,\tilde{v}_1,\ldots, \tilde{v}_{k-1})^T$, $\ul{q} = (q_0,0,\ldots,0)^T$, and $\ulul{A}$ is a $k\times k$ square matrix with $i,j\in\{0,1,\ldots,k-1\}$.  $A$ has $-1$'s along the diagonal and $+1$'s along the first subdiagonal:

\setstretch{1}
\begin{align}
  \ul{\ul{A}} = -\ul{\ul{I}} + \ul{\ul{B}}_k, \qquad
  \ul{\ul{B}}_k = \begin{pmatrix}
    0 &   &        &   &  &  \\
    1 & 0 &        &   &  &  \\
      & 1 &        &   &   &  \\
      &   & \ddots &   &   &  \\
      &   &        & 1 & 0 &  \\
      &   &       &   & 1 & 0 \\
  \end{pmatrix}. \label{eq:Bk}
\end{align}
\setstretch{2}
Note that in successive powers of $\ulul{B}_k$, the $+1$'s move to lower subdiagonals:
\[
  (\ul{B}_k)_{i,j} = \delta_{i-1,j}, (\ul{B}_k^2)_{i,j} = \delta_{i-2,j}, \ldots, (\ul{B}_k^{k-1})_{i,j} = \delta_{i-k+1,j}, \text{ and }\ul{B}_k^l=\ul{0}\text{ for $l\geq k$}.
\]
Equation (\ref{eq:vec}) has the solution
\begin{align}
  \tilde{\ul{v}}(r) = e^{\ulul{A}r}\tilde{\ul{v}}(0)+\int_0^r e^{s\ulul{A}}\ul{q}ds. \label{eq:sol}
\end{align}
Using equation (\ref{eq:Bk}), and the fact that $\ulul{B}_k^{k}=\ulul{0}$, we have that
\[
  e^{Ar} = e^{-r}\ulul{1} e^{-\ulul{B}r}=e^{-r}\left[\ulul{1}+r\ulul{B}_k+\frac{r^2}{2}\ulul{B}^2_k+\cdots +\frac{r^{k-1}}{(k-1)!}\ulul{B}^{k-1}_k\right].
\]
Each term within the square brackets corresponds to a distinct subdiagonal, so the matrix $e^{Ar}$ is 0 above the diagonal and constant along each subdiagonal.  Let $H_j$ be the value in the $j$th subdiagonal: explicitly,
\[
  H_j(r)=e^{-r}\frac{r^j}{j!},
\]
for $j=0,\ldots,k-1$.  Then (\ref{eq:sol}) can be expanded to give an explicit solution for $\tilde{v}_j$ in terms of $r$. 
\begin{align}
  \tilde{v}_j(r) = \sum_{i=0}^j H_{j-i}\tilde{v}_i(0) + \frac{1}{k\lambda}\int_0^r H_j(s)ds, \qquad j=0,\ldots,k-1, \label{eq:explj}
\end{align}
We also have
\begin{align}
  \tilde{v}_k(r) &= 1 - \sum_{j=0}^{k-1} \tilde{v}_j(r) \label{eq:solvk}\\
         &= 1 - \sum_{j=0}^{k-1} \sum_{i=0}^j H_{j-i}(r)\tilde{v}_i(0) - \frac{1}{k\lambda}\int_0^r \sum_{j=0}^{k-1}H_j(s)ds. \nonumber 
\end{align}
Clearly $\tilde{v}_i\rightarrow \bar{v}_i=(\lambda k)^{-1}$ for $i<k$, and $\tilde{v}_k\rightarrow\bar{v}_k = (\lambda -1) /\lambda$ as $r\rightarrow\infty$.  Recalling now that $dr/dt = k\lambda v_k$, we get
\begin{align}
  \frac{dr}{dt} &= k\Big( \lambda - 1 - \lambda  \sum_{j=0}^{k-1}\sum_{i=0}^j H_{j-i}(r)\Big[v_i(0) - \frac{1}{\lambda k}\Big]\Big) = \phi(r), \label{eq:autr}
\end{align}
which is an autonomous ODE, emphasized by introducing the notation $\phi(r)$ on the right.  The behavior of $r(t)$ can then be determined by analyzing $\phi(r)$.  

Starting with $v_k(0)>0$, we see that $\phi(0)>0$, and so $r(t)$ is monotone increasing, and satisfies
\begin{align}
t = \int_0^r \frac{ds}{\phi(s)}
\end{align}
so we have either
\begin{enumerate}
\item there exists $r_0$ the smallest positive solution to $\phi(r)=0$, such that $t\rightarrow\infty$ as $r\rightarrow r_0$, or
\item the integral is finite, in which case $r\rightarrow\infty$ as $t\rightarrow\infty$. This will certainly be the case if $\phi(r)>0$ for all $r>0$.
\end{enumerate}
In case 1, since $\phi(r)=dr/dt=k\lambda v_k$, we see that $v_k \rightarrow 0$ as $t\rightarrow\infty$; that is, in this case the system goes to an inactive state with the firings dying out.  In case 2, we see by equation (\ref{eq:autr}) that as $t\rightarrow\infty$,
\[
v_k(t) \rightarrow \frac{ \phi(\infty) }{ \lambda k } = 1 - \frac{k}{\lambda k} = \bar{v}_k,
\] 
and by equation (\ref{eq:explj}) for $j<k$, that
\[
v_j(t) \rightarrow \tilde{v}_j(\infty) = \frac{1}{k\lambda} = \bar{v}_j,
\]
so the system goes to the unique sustaining stationary state.

\medskip\noindent{\bf \underline{Illustrative Examples}}

\medskip\noindent{\bf The $\mathbf{k=1}$ case.}  In this case, 
\begin{align}
  \frac{dr}{dt}=\lambda v_1(t) = \lambda \left(1-e^{-r}v_0(0)-\frac{1}{\lambda}\int_0^r e^{-s}ds\right)=\phi(r).\label{eq:drdt}
\end{align}
If $0\leq v_0(0) < 1$, then it can be checked that the right-hand side of equation (\ref{eq:drdt}) is bounded below by a positive constant for all $r\geq 0$.  Therefore the system always approaches the sustaining steady state solution if $v_1(0)>0$.  More explicitly, we have in this case, $v_0(t) = 1 - v_1(t)$, so $v_1(t)$ satisfies the autonomous equation
\begin{align}
\frac{dv_1(t)}{dt} = (\lambda -1) v_1(t) -\lambda v_1^2(t),
\end{align}
whose solution is 
\begin{align}
v_1(t) = \frac{ v_1(0) }{\frac{\lambda}{\lambda-1}v_1(0) + \big(1-\frac{\lambda}{\lambda-1}v_1(0)\big)e^{-(\lambda -1)t}} \rightarrow \bar{v}_1 = \frac{\lambda - 1}{\lambda} \text{ as $t\rightarrow \infty$}.
\end{align}
Thus for $k=1$ any initial state with $v_1(0)\neq 0$ will approach, as $t\rightarrow\infty$, the sustaining stationary state exponentially, as long as $\lambda > 1$.  

Thinking of the process as a mean field model of infection, with $v_1$ representing the infected fraction of the population the model predicts as in the standard contact process a persistent percentage of infected individuals for $\lambda > 1$, the percentage increasing with $\lambda$.  For $\lambda < 1$, $v_1(t)\rightarrow 0$, there is no epidemic, as everyone gets eventually cured.

\medskip\noindent{\bf The $\mathbf{k=2}$ case.}  Unlike the $k=1$ case, here it is possible to start the system with $v_2(0)>0$ and still have the firing die out, $v_2(t)\rightarrow 0$ as $t\rightarrow\infty$.

\begin{figure}[h]
\includegraphics[scale=.6]{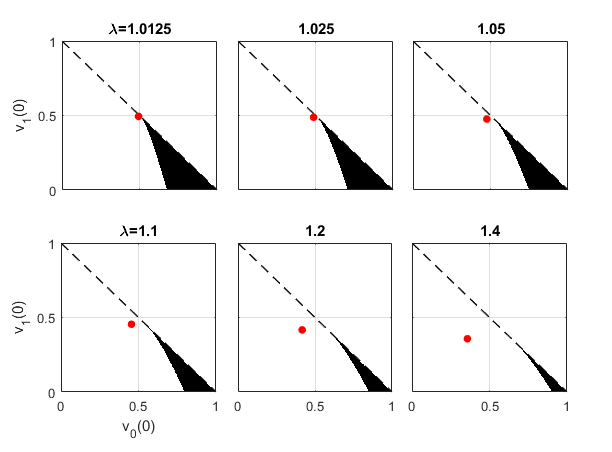}
\caption{\small Solutions of (\ref{eq:29}) for the $k=2$ spatially uniform case, indicating in black the initial conditions $(v_0(0),v_1(0))$ which go to the extinct states $v_2=0$ as $t\rightarrow\infty$.  The extinct states are indicated by the dashed line.  The sustaining stationary state $\bar{v}_0=\bar{v}_1=1/(2\lambda)$ is shown by a red dot: it attracts all initial states in the white region below the dotted line.} \label{fig:3}
\end{figure}

Writing the solution (\ref{eq:explj}) for $k=2$, we get,
\begin{align}
  \tilde{v}_0(r)+\tilde{v}_1(r) &=\frac{1}{\lambda}[1-e^{-r}]-\frac{1}{2\lambda}re^{-r}+v_0(0)e^{-r}[1+r] + v_1(0)e^{-r}. \nonumber
\end{align}
This yields
\begin{align}
\phi(r) &= 2\lambda [1 - \tilde{v}_1(r) - \tilde{v}_0(r)] \label{eq:29} \\
        &= 2(\lambda -1)(1-e^{-r}) + re^{-r} [1 - 2\lambda v_0(0)] + 2\lambda e^{-r}[1-v_1(0)-v_0(0)] \nonumber. 
\end{align}
Clearly if $v_0(0)<\frac{1}{2\lambda}$ then $\phi(r)>0$ for all $r>0$ and the system will go to the sustaining stationary state.  For $v_0(0)$ close to 1, there exists $r_0>0$ for which $\phi(r_0)=0$, in which case the system goes to an inactive state with the firing dying out, i.e., $v_2\rightarrow 0$ as $t\rightarrow \infty$.  The region in the $v_0(0)$, $v_1(0)$ plane for which this occurs shrinks as $\lambda$ increases, as demonstrated numerically in figure \ref{fig:3}.  We expect similar behavior for $k>2$.  In fact it is easy to see from (\ref{eq:autr}) that if $v_i(0)\leq \frac{1}{\lambda k}$ for all $i\leq k-1$ then $\phi(r) > 0$.

\section{Linear Stability of the Sustaining State}\label{sec:5}
Let us consider now the linear stability of the sustaining stationary states for general $k$
\begin{align*}
  \bar{v}_j = \frac{1}{\lambda k}, \qquad j<k, \qquad
  \bar{v}_k = 1-\frac{1}{\lambda}.
\end{align*}
Let $v_j(\ul{x},0)=\bar{v}_j+f_j(\ul{x},t)$.  Then setting
\[
  v_j(\ul{x},t)=\bar{v}_j + f_j(\ul{x},t), 
\]
The conditions 
\begin{align*}
&\sum_{j=0}^k f_j(\ul{x},t)=0, \\
&1-\bar{v}_j \geq f_j(\ul{x},t) \geq -\bar{v}_j
\end{align*}
are satisfied for all $t\geq 0$ if they are satisfied for $t=0$.

Setting $(\lambda -1)k=\alpha$, we obtain to first order in the $f_j$, 
\begin{align}
  \frac{\partial f_0}{\partial t}(\ul{x},t) &= -\alpha f_0(\ul{x},t) + \int J(\ul{x}-\ul{y})[f_k(\ul{x})-f_k(\ul{y})]d\ul{y} \label{eq:lineq1}\\
  \frac{\partial f_j}{\partial t}(\ul{x},t) &= \alpha [f_{j-1}(\ul{x},t)-f_j(\ul{x},t)], \qquad j=1,\ldots,k-1. \label{eq:lineq2}
\end{align}
Taking the Fourier series in $\ul{x}$,
\begin{align}
\hat{f}_j(\ul{\xi},t) = \frac{1}{L^d}\int_\Lambda f_j(\ul{x},t)e^{-2\pi i \ul{\xi} \cdot \ul{x}} d\ul{x}, \qquad j=0,\ldots,k-1,
\end{align}
where $\ul{\xi}=(\xi_1,\ldots,\xi_d) \in (L^{-1}\mathbb{Z})^d$, gives
\begin{align}
\frac{\partial \hat{f}_0}{\partial t} (\ul{\xi},t) &= -\alpha \hat{f_0}(\ul{\xi},t) + \hat{f}_k (\ul{\xi},t) \big(1 - \hat{J}(\ul{\xi})\big) \\
\frac{\partial \hat{f}_j}{\partial t} (\ul{\xi},t) &= \alpha [\hat{f}_{j-1}(\ul{\xi},t) - \hat{f}_j(\ul{\xi},t)], \qquad j=1,\ldots, k-1.
\end{align}
Letting $\beta(\ul{\xi}) = (1-\hat{J}(\ul{\xi})) / \alpha$, we can write the above equations in vector form
\begin{align}
\frac{\partial \hat{\ul{f}}}{\partial t} (\ul{\xi},t) = \alpha (\ul{A} - \beta(\ul{\xi})\ul{M})\hat{\ul{f}}, \qquad \ul{M} =  
\begin{pmatrix}
  1  & 1 & \cdots & 1 \\
  0  & 0 & \cdots & 0 \\
  \vdots & & \ddots & \vdots \\
  0  & 0 & \cdots & 0 
\end{pmatrix}, \label{eq:linsys}
\end{align}
where $\ul{M}$ is the $k \times k$ matrix with ones in the top row and zeros elsewhere.

We will show that for any fixed value of $\ul{\xi}$, each eigenvalue of the matrix $(\ul{A}-\beta\ul{M})$ has a negative real part, implying the convergence of $\hat{\ul{f}}(\ul{\xi},t)$ to $\ul{0}$ as $t\rightarrow\infty$.  The characteristic polynomial of the matrix, computed by cofactor expansion along the top row, is
\begin{align}
p(x) &= (-1 -\beta - x) - \sum_{k=2}^k (-1)^j(-\beta)(+1)^{j-1}(-1-x)^{k-j} \\
     &= (-1)^k\Big[(1+\beta+x)(1+x)^{k-1}+\beta\sum_{j=2}^k(1+x)^{k-j}\Big]. \nonumber
\end{align}
Letting $p(x)=0$, $y=x+1$, and noting that $y= 1$ is not a solution to the equation, yields
\begin{align}
0 &= (y+b)y^{k-1}+\beta\sum_{j=2}^k y^{k-j} \\
 &= y^k + \beta\sum_{j=1}^k y^{k-j} \nonumber \\
&= y^k + \beta\frac{1-y^k}{1-y}, \nonumber
\end{align}
which simplifies to
\begin{align}
y^k(y-(1-\beta)) = \beta. \label{eq:charsimp}
\end{align}
We now note that $\beta\geq 0$. This follows from the fact that $J$ is a real-valued function with $J(\ul{x})=J(-\ul{x})\geq 0$ and $\int_\Lambda J(\ul{x})d\ul{x}=1$.  Therefore applying the absolute value to both sides of (\ref{eq:charsimp}) gives
\begin{align}
|y|^k|y-(1-\beta)| = \beta. \label{eq:charsimp2}
\end{align}
We can deduce that $|y|\neq 1$, as otherwise the equation (\ref{eq:charsimp2}) implies $y=1$, and we have remarked that this not a solution to the characteristic equation.  Therefore either 
\begin{enumerate}
\item $|y|<1$, or
\item $|y|>1$.
\end{enumerate}
In case 1, it follows easily that the real part of $y$ is less than 1, and therefore the real part of $x$ is negative.  Case 2 implies that $|y-(1-\beta)|<\beta$, from which it also follows that the real part of $y$ is less than 1, and that the real part of $x$ is negative.  This implies that $\hat{f}_j(\ul{\xi},t)\rightarrow 0$ for all fixed $\xi$, and from this point, we can prove that $\ul{f}$ converges to $\ul{0}$ uniformly.

\medskip\noindent{\bf Linear stability of the inert state $\mathbf{v_0=1}$.}  
Having shown the linear stability of the stationary sustaining state we consider now the linear stability of the state $v_k(t) = v_k(0)=0$. Rather than considering all initial states $v_j(0)$, $\sum_{j=0}^{k-1}v_j(0)=1$, we consider here only perturbations around the extreme case $v_0(\ul{x},t)=1$, $v_j(\ul{x},t)=0$, $j>0$. 

Let
\begin{align}
v_0(\ul{x},t) &= 1 + f_0(\ul{x},t), \qquad f_0(\ul{x},t) < 0, \\
v_j(\ul{x},t) &= f_j(\ul{x},t), \qquad f_j(\ul{x},t) > 0, \qquad j=1,\ldots,k, \qquad \sum_{j=0}^k f_j=0.
\end{align}
Linearizing in the $f$'s gives, for $k>1$,
\begin{align}
\frac{\partial f_0}{\partial t} (\ul{x},t) &= f_k(\ul{x},t) - \lambda k\int J(\ul{x}-\ul{y})f_k(\ul{y},t)d\ul{y} \\
\frac{\partial f_1}{\partial t} (\ul{x},t) &= \lambda k\int J(\ul{x}-\ul{y})f_k(\ul{y},t)d\ul{y} \\
\frac{\partial f_j}{\partial t} (\ul{x},t) &= 0, \quad \text{for $j\in\{2,\ldots,k-1\}$}, \\
\frac{\partial f_k}{\partial t} (\ul{x},t) &= -f_k = \sum_{j=0}^{k-1} f_j, \quad \text{ implying } f_k(\ul{x},t)=f_k(\ul{x},0)e^{-t},
\end{align}
so for $k>1$ the dead state in the vicinity of $v_0(\ul{x},t)=1$ is linearly stable, see figure \ref{fig:3}.

For $k=1$, $f_1(\ul{x},t)\geq 0$, 
\begin{align}
\frac{\partial f_1}{\partial t} (\ul{x},t) = -f_1(\ul{x},t) + \lambda \int J(\ul{x}-\ul{y})f_1(\ul{y},t)d\ul{y}.
\end{align}
Taking spatial Fourier transforms yields
\begin{align}
\frac{\partial \hat{f}_1}{\partial t} (\ul{\xi},t) &= (\lambda \hat{J}(\ul{\xi}) -1) \hat{f}_1(\ul{\xi},t),
\end{align}
which has the solution
\begin{align}
\hat{f}_1(\ul{\xi},t) &= \hat{f}_1(\ul{\xi},0)e^{(\lambda \hat{J}(\ul{\xi}) -1)t}.
\end{align}
Since $\hat{J}(0)=1$ there will be growth at least for small values of $\ul{\xi}$, for which $\hat{f}_1(\ul{\xi},0) > 0$, so the state $v_0=1$ is unstable for $\lambda > 1$.  In fact as we shall now show for $k=1$, any perturbation of the state $v_0(0)=1$, will lead asymptotically to the stable stationary state $\bar{v}_0=\lambda^{-1}$, $\bar{v}_1 = (\lambda -1)/\lambda.$

\section{$k=1$, general case.}\label{sec:6} In this two-level case, $v_1(\ul{x},t)=1-v_0(\ul{x},t)$ is the only unknown function.  It satisfies the equation 
\begin{align}
\frac{\partial v_1(\ul{x},t)}{\partial t} = -v_1 + \lambda (1-v_1(\ul{x},t)) \int J(\ul{x}-\ul{y}) v_1(\ul{y},t)d\ul{y}. \label{eq:3.1}
\end{align}
Define $f(\ul{x},t)$ by
\begin{align}
f(\ul{x},t) = v_1(\ul{x},t) - \bar{v}_1, \qquad \bar{v}_1 = \frac{\lambda -1}{\lambda}.
\end{align}
Using the fact that $\int J(\ul{x}-\ul{y})d\ul{y}=\int J(\ul{x}-\ul{y})d\ul{x}=1$ we get from (\ref{eq:3.1})
\begin{align}
\frac{\partial f(\ul{x},t)}{\partial t} = &-(\lambda -1)f(\ul{x}) \label{eq:3.6} \\ &+\int_\Lambda J(\ul{x}-\ul{y})[f(\ul{y})-f(\ul{x})]d\ul{y} \nonumber \\&- \lambda \int_\Lambda J(\ul{x}-\ul{y}) f(\ul{y})f(\ul{x})d\ul{y}. \nonumber \\
&= -\int d\ul{y} J(\ul{x}-\ul{y})[f(\ul{x})-f(\ul{y})] - \lambda\int d\ul{y}J(\ul{x}-\ul{y})f(x)\left[ \frac{\lambda -1}{\lambda}+f(\ul{y})\right] \nonumber
\end{align}
Multiplying (\ref{eq:3.6}) by $f(\ul{x},t)$ and integrating over $\ul{x}$ yields
\begin{align}
\frac{1}{2}\frac{d}{dt} \int_\Lambda f^2(\ul{x},t) d\ul{x} = &-\frac{1}{2} \iint J(\ul{x}-\ul{y})[f(\ul{x})-f(\ul{y})]^2 d\ul{x}d\ul{y}  \label{eq:25} \\ &- \lambda \iint J(\ul{x}-\ul{y}) f^2(\ul{x},t) v_1(\ul{y},t) d\ul{y}d\ul{x} \leq 0 \nonumber
\end{align}
The inequality is strict for all initial conditions with $v_1(\ul{x},0)$ not identically 0, and shows that $v_1(\ul{x},t)\rightarrow \bar{v}_1$ as long as $v_1(\ul{x},0)>0$.

\begin{figure}[h]
\includegraphics[scale=.14]{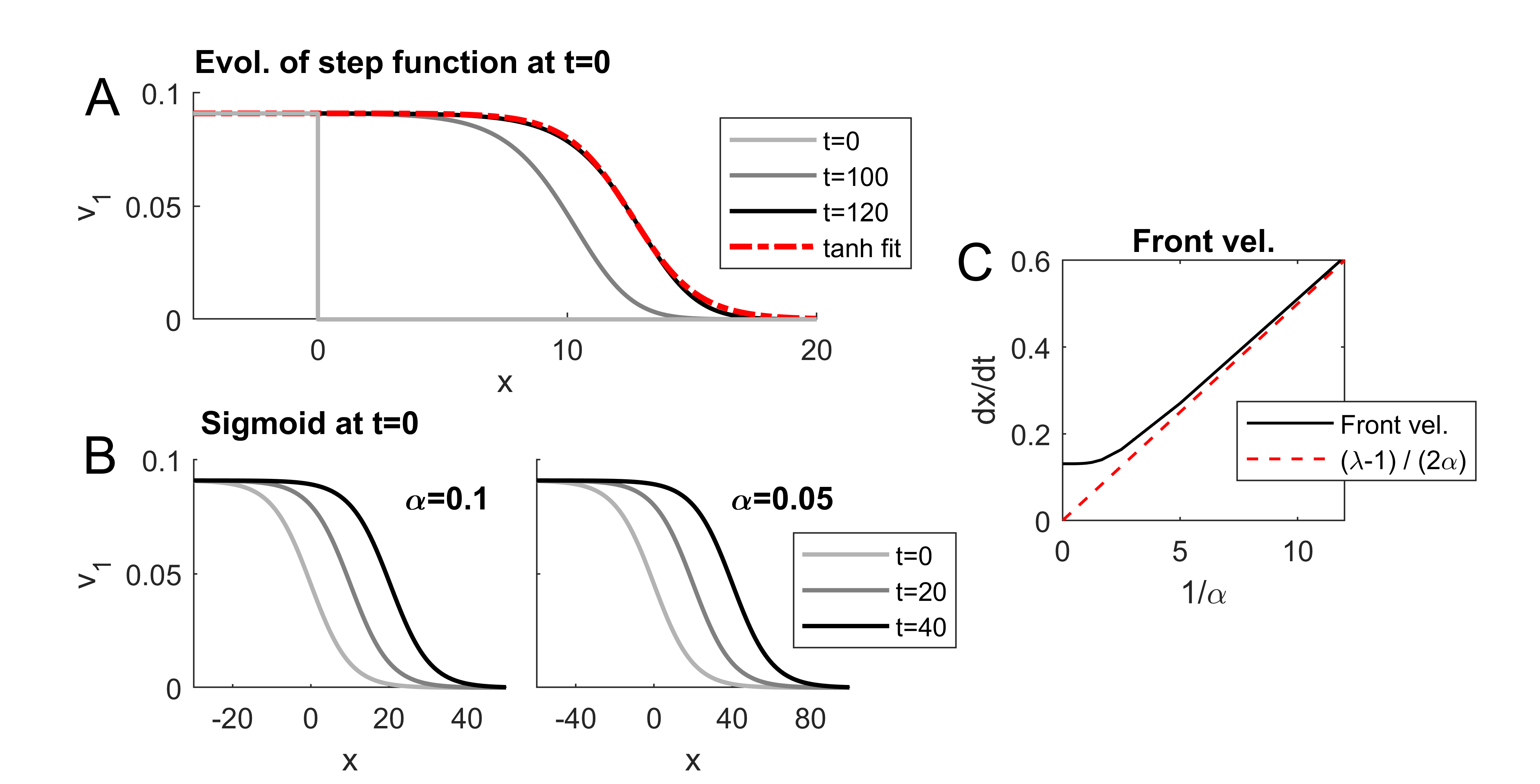}
\caption{\small Numerical simulation of $k=1$ traveling front solutions.  In all plots, $\lambda=1.1$ and $J(x)=1$ for $|x|<1/2$ and 0 otherwise.  {\bf A.} Traveling front arising from step function initial conditions. At $t=0$, $v_1(x)=\bar{v}_1$ for $x<0$ and 0 otherwise. The traveling front is well approximated by a rescaled and shifted hyperbolic tangent function, shown in red. {\bf B.} Fronts developing from initial conditions set to eq. (\ref{eq:trav}) with $\alpha=0.1$ and $\alpha=0.05$.  Front velocity is approx.\ double for $\alpha=0.05$ compared to $0.1$.  {\bf C.} Dependence of front velocity on $\alpha$. $v_1(x)$ is initialized with eq. (\ref{eq:trav}) for a range of $\alpha$, and velocity is computed numerically.  Comparison with velocity given in eq. (\ref{eq:trav}) is given by the red dashed line.  Note that for $1/\alpha$ close to 0, the initial condition is nearly a step function as in panel {\bf A}. }\label{fig:1}
\end{figure}

\medskip\noindent{\bf Traveling wave solution.}
We consider the time evolution of $v_1(x,t)$ when $x\in\mathbb{R}$, i.e., we let $L\rightarrow\infty$, and the initial state is one in which $v_1(x,0)$ goes to the stable solution $\bar{v}_1=(\lambda-1)/\lambda$ as $x\rightarrow -\infty$ and to the unstable solution $v_1=0$ as $x\rightarrow\infty$.  Equation (\ref{eq:3.1}) for $k=1$ can be considered as a special case of the non-local KPP equation, with the diffusion constant set equal to zero \cite{Berestycki}.  To get a feeling for the evolution of such an initial state we first consider the limiting case when the width $J(x-y)$ goes to 0, i.e.,
\begin{align}
J(x-y)=\delta (x-y).
\end{align}
Equation (\ref{eq:3.1}) then has the traveling wave solution 
\begin{align}
u_1(x,t) = \frac{\lambda -1}{\lambda} \Big[1-\tanh [\alpha (x-Vt)]\Big]/2, \quad \alpha V = (\lambda -1)/2 \label{eq:trav}
\end{align}
Numerical solutions of (\ref{eq:3.1}) with 
\begin{align*}
J(x)=\frac{1}{2b}\theta(b-|x|), \quad \theta(x)=\begin{cases}1, &\quad x\geq0, \\0, &\quad x<0,\end{cases}
\end{align*}
and initial conditions 
\[
v_1(x,0)=\begin{cases}
\frac{\lambda-1}{\lambda}, \quad &\text{for $x < 0$},\\
0, \quad &\text{for $x > 0$},
\end{cases}
\]
show that $v_1(x,t)$ approaches a form close to (\ref{eq:trav}) with $\alpha V\sim (1-\lambda)/2$ as $t\rightarrow \infty$; see figure \ref{fig:1}.  Similar behavior is found for $J(x)$ a Gaussian.

\section{Conclusions}\label{sec:concl}

We have shown that when $\lambda > 1$, the macroscopic equations (\ref{eq:4})-(\ref{eq:6}) have a unique, linearly-stable stationary state with nonzero firing rates $v_k(\ul{x},t)$ given by (\ref{eq:13p}), which we refer to as the sustaining stationary state.  For $k=1$, the basin of attraction of the sustaining stationary state includes all initial conditions such that $v_1(\ul{x},0)$ is not identically $0$ for all $x$.  The case $k>1$ is qualitatively different: even starting with $v_k(0)>0$ we can have $v_k(t)\rightarrow 0$.  We have shown this explicitly for $k=2$ and found that there are linearly-stable extinct states like $\ul{v}(x,t)=(1,0,\ldots,0)$.  

For $k=2$, we have seen in numerical simulations that even for initial conditions with $\ul{v}(x,0)=\bar{\ul{v}}$ for $x$ in a small region of $\mathbb{R}$ and $\ul{v}(x,0)=(1,0,0)$ outside of that region ($J$ the same as in figure~\ref{fig:1}), the firing (epidemic) can die out and approach an extinct state.  As the size of the region increases, eventually a point is reached where the firing becomes self-sustaining and spreads throughout the system.  We conjecture that there exists some $M>0$, depending on $J(x)$, such that the initial condition with the region $|\ul{x}|<M$ set to $\bar{\ul{v}}$ will necessarily converge to $\bar{\ul{v}}$ pointwise on the whole domain.

In the one-dimensional case $x\in\mathbb{R}$, we showed the existence of traveling wavefront solutions $\ul{v}(x,t)=\ul{v}(x-Vt)$ in the $k=1$ case, with an analytic solution in the $J=\delta$ case and numerically for other forms of $J$.  Wavefront solutions have been studied extensively in the case of neural models like the mean field Wilson-Cowan equations \cite{Ramirez}.  It remains to be rigorously shown that traveling fronts exist in the case $k>1$.  Based on numerical solutions of the equations, we conjecture that stable traveling wave solutions exist for $k>1$, and that there are qualitative differences with the case $k=1$.  The existence of linearly stable extinct states like $\ul{v}=(1,0\ldots,0)$ when $k>1$ changes the properties of traveling fronts with $\ul{v}(-\infty)=\bar{\ul{v}}$ and $\ul{v}(\infty)=(1,0,\ldots,0)$.  In particular, we conjecture that there exists unique traveling waves, with wave velocity proportional to the width of $J$.  

As the microscopic system described in the Introduction has been inspired by integrate and fire models from neuroscience, we are interested in extensions of the model to include more realistic components of biological neurons and neural networks.  In particular, including sites in the model which have an inhibitory effect on nearby sites will lead to a richer dynamical landscape.  Other models with mixed excitatory and inhibitory components have been shown to have oscillatory activity, and such models have been used to study brain rhythms (e.g., \cite{Brunel}, \cite{Keeley}, \cite{Wallace}).  Additionally, the effects of inhibition on the propagation of traveling fronts in the brain has been examined in \cite{Ramirez}.  An interesting question is how brain activity is contained and localized to a particular region when externally driven.  We would like to explore this question in our model with the addition of both inhibition and an external drive component.










\medskip\noindent{\bf Acknowledgements.} We would like to thank Lai-Sang Young, Anna De Masi, and Errico Presutti for many helpful discussions.

\end{document}